\documentclass[12pt,a4paper]{article}
\usepackage{pifont}
\usepackage{amsfonts}
\usepackage{mathrsfs}
\usepackage{amssymb}
\usepackage{amsmath}

\newtheorem{theorem}{Theorem}[section]
\newtheorem{definition}[theorem]{Definition}
\newtheorem{lemma}[theorem]{Lemma}

\makeatletter

\newcommand{\Rmnum}[1]{\expandafter\@slowromancap\romannumeral #1@}
\makeatother

\begin{document}
\title{Gr\"{o}bner-Shirshov bases for free partially commutative  Lie algebras
\footnote{Supported by the NNSF of China (Nos. 10771077;
10911120389).}}
\author{
 Yuqun
Chen \ and   Qiuhui Mo \\
{\small \ School of Mathematical Sciences, South China Normal
University}\\
{\small Guangzhou 510631, P.R. China}\\
{\small yqchen@scnu.edu.cn}\\
{\small scnuhuashimomo@126.com}}

\date{}

\maketitle \noindent\textbf{Abstract:}  In this paper, by using
Composition-Diamond lemma for Lie algebras, we give a
Gr\"{o}bner-Shirshov basis for free partially commutative Lie
algebra over a commutative ring with unit. As an application, we
obtain a normal form for such a Lie algebra.



\section{Introduction}
The free partially commutative monoid was introduced by P. Cartier
and D. Foata in 1969 \cite{CaFo} for the study of combinatorial
problems in connection with word rearrangements. Since that time,
this monoid has been the subject of many studies. They were
principally motivated by the fact that the free partially
commutative monoid is a model for concurrent computing. On the other
hand, they can be seen as a natural generalization of free monoids.
Indeed, several classical results of the free monoid theory can be
extended to the partially commutative framework. For instance, the
free partially commutative group, the free partially commutative
associative algebra and the free partially commutative Lie algebra.
This research direction was followed by a lot of people (see
\cite{Cho, Di, Db.1, Du.1,  Duch, Sm, Th} or \cite{Pe}).

In 1992, G. Duchamp \cite{Duch} proved that the free partially
commutative Lie algebra is a free $K$-module, where $K$ is a
commutative ring, and gave an algorithm to find a basis for any free
partially commutative Lie algebras. In 1993, G. Duchamp and D. Krob
\cite{DuKr} showed how to obtain decomposition results for free
partially commutative Lie algebra into free Lie algebras and
obtained a normal form for such algebra. The algorithms given by
these two papers are based on a decomposition of the generating set
of the  free partially commutative Lie algebra by two subsets one of
which is independent, and a linear basis is implicitly given.

In 2001, by using Composition-Diamond lemma for associative
algebras, L.A. Bokut and L.S. Shiao  \cite{BoSh} gave
Gr\"{o}bner-Shirshov bases for free partially commutative
associative algebras and free partially commutative monoids. In
2004, E.S. Esyp, I.V. Kazatchkov and V.N. Remeslennikov \cite{Esyp}
gave a Gr\"{o}bner-Shirshov basis for free partially commutative
group. In this paper, by using Composition-Diamond lemma for Lie
algebras, we give a Gr\"{o}bner-Shirshov basis for free partially
commutative Lie algebra over a commutative ring with unit. As an
application, we obtain a linear basis explicitly for such a Lie
algebra. E.N. Poroshenko \cite{poro} also gives the same
Gr\"{o}bner-Shirshov basis for free partially commutative Lie
algebra independently.  Our proof is different from \cite{poro}.

Let $X$ be a set, $K$ a commutative ring with unit and $Lie(X)$ the
free Lie algebra over $K$ generated by $X$. Let
$\vartheta\subseteq(X\times X)\backslash \{(x,x) \mid x\in X\}$.
Then $Lie(X|\vartheta)=Lie(X)/Id(\vartheta)$ is free partially
commutative Lie algebra, where $Id(\vartheta)$ is the ideal of
$Lie(X)$ generated by the set $\{(ab)| (a,b)\in\vartheta\}$ and
$(ab)$ is the Lie multiplication in $Lie(X)$.

\section{Preliminaries}
We start with the associative Lyndon-Shirshov words.

Let $X=\{x_i|i\in I\}$ be a well-ordered set with $x_i>x_p$ if $i>p$
for any $i,p\in I$ and $X^*$ the free monoid generated by $X$.  We
order $X^*$ by the lexicographical ordering.

\begin{definition}(\cite{BoC1,Chen,Ly,Sh62,Shir3,Uf})
An associative word  $u\ (u\neq 1)$ in $X$ is called an $ALSW$
(associative Lyndon-Shirshov word) if
$$
(\forall v,w\in X^\ast, v,w\neq 1) \ u=vw\Rightarrow vw>wv.
$$

A non-associative word  $(u)$ in $X$ is called a $NLSW$
(non-associative Lyndon-Shirshov word)  if
\begin{enumerate}
\item[(i)] $u$ is an $ALSW$,
\item[(ii)] if $(u)=((v)(w))$, then both $(v)$ and $(w)$ are
$NLSW$'s,
\item[(iii)] in (ii) if $(v)=((v_1)(v_2))$, then $v_2\leq w$ in
$X^\ast$.
\end{enumerate}
\end{definition}

\begin{lemma}\label{l2.15}(\cite{Chen,Sh62,Shir3})
(i) For any $ u \in X^*$, there exists a unique decomposition
$u=u_1u_2 \cdots u_k$, where $u_i$ is an $ALSW$, \ $1 \leq i \leq
k$, \ \ and \ $u_1 \leq u_2 \leq \cdots \leq u_k$.

(ii) Let $u$ be an $ALSW$ \ and $|u| \geq 2$. If $u=vw$, where $w$
is the longest $ALSW$ proper end of $u$, then $v$ is an $ALSW$.
\end{lemma}

Following \cite{Chen,Sh62}, for an $ALSW$ $u$, there is a unique
bracketing, denoted by $[u]$, such that $[u]$ is $NLSW$:
$$
[x_i]=x_i, \ [u]=[[v][w]],
$$
where $u=vw$ and $w$ is the longest $ALSW$  proper end of $u$.

\ \

Now, we consider $(\ )$ as Lie bracket in the free associative
algebra $k\langle X\rangle$, i.e., for any $a,b\in k\langle
X\rangle, \ (ab)=ab-ba$, where $k$ is a field. We may view $Lie(X)$
as the subLie-algebra of $k\langle X\rangle$ generated by $X$.

For any polynomial $f\in k\langle X\rangle$, $f$ has the leading
word $\overline{f}$. We call $f$  monic if the coefficient of
$\overline{f}$ is 1. By $deg(f)$ we denote the degree of
$\overline{f}$.

\begin{lemma}(\cite{BoC1,Chen,Ly,Sh62,Shir3,Uf})\label{l3.6}
$NLSW$'s forms a linear basis of $Lie(X)$.
\end{lemma}

\begin{lemma}(\cite{Chen,Sh62,Shir3})\label{l3.5}
Let $[u]$ be a NLSW. If we consider $[u]$ as a polynomial in
$k\langle X\rangle$,  then $\overline{[u]}=u$.
\end{lemma}

\begin{lemma}(\cite{Sh62,Shir3})\label{l3.0}
Let $u,v$ be $ALSW$'s, $u=avb, \ a,b\in{X^*}$. Then $[u]=[a[vc]d]$,
where $b=cd, \ c,d\in{X^*}$. Denote by
$$
[u]_v=[u]|_{[vc]\mapsto{[[[v][c_1]]\cdots[c_k]]}},
$$
where $c=c_1 \cdots c_k$, $c_j$ is an $ALSW$ and $c_1 \leq{c_2
}\leq\cdots \leq{c_k}$. Then $\overline{[u]_v}=u$.
\end{lemma}

Let $S\subset Lie(X)$ with each $s\in S$ monic, $a,b\in{X^*}$ and
$s\in S$. If $a\bar{s}b$ is an ALSW, then we define
$[asb]_{\bar{s}}=[a\bar{s}b]_{\bar{s}}|_{[\bar{s}]\mapsto{s}}$,
where $[a\bar{s}b]_{\bar{s}}$ is defined by Lemma \ref{l3.0} (see
\cite{BoC1}).

From now on, we use the deg-lex ordering $<$ on $X^*$: to compare
two words by degree first and then lexicographically.

Let $f$ and $g$ be two monic Lie polynomials in $Lie(X)\subset
k\langle X\rangle$. Then, there are two kinds of Lie compositions:
\begin{enumerate}
\item[(i)] If  $w=\bar{f}=a\bar{g}b$ for some $a,b\in X^*$, then the
polynomial $(f,g)_w=f - [agb]_{\bar{g}}$ is called the
 composition of inclusion of $f$ and $g$ with respect to $w$.

\item[(ii)] If \ $w$ is a word such that $w=\bar{f}b=a\bar{g}$ for
some $a,b\in X^*$ with $deg$$(\bar{f})+$deg$(\bar{g})>$deg$(w)$,
then the polynomial
 $(f,g)_w=[fb]_{\bar{f}}-[ag]_{\bar{g}}$ is called the composition of intersection of $f$ and
$g$ with respect to $w$.

The $w$ in the above is called ambiguity.
\end{enumerate}

Let $S\subset Lie(X)$ with each $s\in S$ monic.

Suppose that $a,b\in{X^*}$ and $s\in S$. If
$\overline{(asb)}=a\overline{s}b$ and $a\bar{s}b$ is an $ALSW$, then
we call $(asb)$ a normal $s$-word (or normal $S$-word).

Suppose that $w\in X^*$ and $h$ is a Lie polynomial. Then  $h$ is
trivial modulo $(S,w)$, denoted by $h\equiv0 \ mod(S,w)$, if
$h=\sum\limits_{i}\alpha_i(a_is_ib_i)$, where each
$\alpha_{i}\in{k}$, $a_i, b_i\in{X^*}$, $s_i\in{S}$, and
$(a_is_ib_i)$ is normal $S$-word such that $a_i\overline{s_i}b_i<w$.

The set $S$ is called a Gr\"{o}bner-Shirshov basis in $Lie(X)$ if
any composition in $S$ is trivial modulo $S$ and corresponding $w$.

\begin{lemma}(\cite{Sh62,Shir3}, Composition-Diamond lemma for Lie algebras)\label{cdL}
Let $S\subset{Lie(X)}$ be nonempty set of monic Lie polynomials. Let
$Id(S)$ be the ideal of $Lie(X)$ generated by $S$. Then the
following statements are equivalent.
\begin{enumerate}
\item[(i)] $S$ is a Gr\"{o}bner-Shirshov basis in
$Lie(X)$.
\item[(ii)] $f\in{Id(S)}\Rightarrow{\bar{f}=a\bar{s}b}$ for
some $s\in{S}$ and $a,b\in{X^*}$.
\item[(iii)]$Irr(S)=\{[u] \ | \ [u] \mbox{ is a } NLSW, \ u\neq{a\bar{s}b}, \
s\in{S},\ a,b\in{X^*}\}$ is a $k$-basis for $Lie(X|S)=Lie(X)/Id(S)$.
\end{enumerate}
\end{lemma}

\noindent{\bf Remark}: In this section, if the field $k$ is replaced
by a commutative ring $K$ with unit, then all results hold. In
particular, the Lemma \ref{cdL} is true for the free Lie algebra
over $K$.

\section{Free partially commutative Lie algebra}

Let  $<$ be a well ordering on $X$.  Throughout this paper,  if
$a>b$ and $(a,b)\in \vartheta$ or $(b,a)\in \vartheta$, we denote
$a\rhd b$. Generally, for any set $Y$, $a\rhd Y$ means $a\rhd y$ for
any $y\in Y$ and $a> Y$ means $a> y$ for any $y\in Y$. For any
$u=x_{i_1}\cdots x_{i_n}\in X^*$ where $x_{i_j}\in X$, we denote the
set $\{x_{i_j}, j=1,\ldots, n\}$ by $supp(u)$.

For any $u\in X^*$, we introduce the two following notions of
degree:
\begin{enumerate}
\item[---] the partial degree $|u|_x$ of $u$ in $x\in X$ is just the
number of $x$ in $u$.

\item[---] the multidegree  is  the $X$-uple $|u|_X=(|u|_x)_{x\in X}\in
N^{(X)}$, where $N$ is the set of  non-negative integers.
\end{enumerate}

\begin{lemma}\label{5.1}
Let $u\in X^*$ be an $ALSW$ and $x\in X$ such that $x> supp(u)$.
Then, in $Lie(X)$, $[xu]=[x[u]]=\sum _{i=1}^m \alpha_i((xy_i)u_i)$,
where for any $i=1,\ldots,m$, $\alpha_i\in k, y_i\in X, u_i\in X^*,
supp(xu)=supp(xy_iu_i)$, $|xu|_X=|xy_iu_i|_X$ and
$\overline{((xy_i)u_i)}=xy_iu_i$.
\end{lemma}

\textbf{Proof}\ We prove the lemma by induction on $|u|$. If
$|u|=1$, the result is clear. Let us suppose it has been proved for
$|u|<n$ with $n\geq 2$. Let $|u|=n$, and $u=u_1u_2$ where $u_2$ is
the longest $ALSW$ proper end of $u$. Then by Lemma \ref{l2.15},
$u_1$ is an $ALSW$ and
$[xu]=(x([u_1][u_2]))=((x[u_1])[u_2])-((x[u_2])[u_1])$. By
induction, $(x[u_1])=\sum _{i=1}^l \beta_i((xy_i)v_i)$, where for
any $i=1,\ldots,l$, $\beta_i\in k, y_i\in X, v_i\in X^*,
supp(xu_1)=supp(xy_iv_i)$, $|xu_1|_X=|xy_iv_i|_X$,
$\overline{((xy_i)v_i)}=xy_iv_i$ and $(x[u_2])=\sum _{j=1}^t
\gamma_j((xz_j)w_j)$, where for any $j=1,\ldots,t$, $\gamma_j\in k,
z_j\in X, w_j\in X^*, supp(xu_2)=supp(xz_jw_j)$,
$|xu_2|_X=|xz_jw_j|_X$,  $\overline{((xz_j)w_j)}=xz_jw_j$. Then
$[xu]=((x[u_1])[u_2])-((x[u_2])[u_1])=\sum _{i=1}^l
\beta_i(((xy_i)v_i)[u_2])-\sum _{j=1}^t \gamma_j(((xz_j)w_j)[u_1])$
and $supp(xy_iv_iu_2)=supp(xu_1u_2)=supp(xu)$,
$supp(xz_jw_ju_1)=supp(xu_2u_1)=supp(xu)$,
$|xy_iv_iu_2|_X=|xu_1u_2|_X=|xu|_X$,
$|xz_jw_ju_1|_X=|xu_2u_1|_X=|xu|_X$. Since $x> supp(u)$ and $u_1,
u_2$ are $ALSW$,  we get that
$\overline{(((xy_i)v_i)[u_2])}=\overline{((xy_i)v_i)}\cdot\overline{[u_2]}=xy_iv_iu_2$
and
$\overline{(((xz_j)w_j)[u_1])}=\overline{((xz_j)w_j)}\cdot\overline{[u_1]}=xz_jw_ju_1$.
\hfill $\blacksquare$

\begin{lemma}\label{5.14}
Let  $x,y,z\in X$ and $u,v\in X^*$ such that $x> y > supp(u)$ and
$y> z
>supp(v)$. Then the following statements hold.
\begin{enumerate}
\item [(i)]\ If $v$ is an $ALSW$, then
$(([xuy][v])z)-([xu]((y[v])z))=(([xu][v])(yz))+((([xu]z)y)[v])-(([xu](z[v]))y)$.
\item [(ii)]\ Suppose that $v=v_1v_2 \cdots v_n$, where $n\geq2$, $v_i$
is an $ALSW$, \ $1 \leq i \leq n$\ and \ $v_1 \leq v_2 \leq \cdots
\leq v_n$. Then
\begin{eqnarray*}
&&((((([xu]y)[v_1])[v_2]) \cdots [v_n])z)-([xu]((((y[v_1])[v_2])
\cdots [v_n])z))\\
&=&(((((([xu]y)[v_1])[v_2]) \cdots
[v_{n-1}])z)[v_n])-(([xu]((((y[v_1])[v_2]) \cdots
[v_{n-1}])z))[v_n])\\
&&+(([xu][v_n])((((y[v_1])[v_2]) \cdots
[v_{n-1}])z))-(([xu](z[v_n]))(((y[v_1])[v_2])
\cdots [v_{n-1}]))\\
&&-\sum_{i=1}^{n-1}((((([xu][v_i])((((y[v_1])[v_2]) \cdots
[v_{i-1}])))[v_{i+1}])\cdots[v_{n-1}])(z[v_n])).
\end{eqnarray*}
\end{enumerate}
\end{lemma}
\textbf{Proof}\ (i)
\begin{eqnarray*}
&&(([xuy][v])z)-([xu]((y[v])z))\\
&=&((([xu]y)[v])z)-(([xu](y[v]))z)+(([xu]z)(y[v]))\\
&=&((([xu]y)[v])z)-((([xu]y)[v])z)+((([xu][v])y)z)+((([xu]z)y)[v])-((([xu]z)[v])y)\\
&=& ((([xu][v])z)y)+(([xu][v])(yz))+((([xu]z)y)[v])-((([xu][v])z)y)-(([xu](z[v]))y)\\
&=&(([xu][v])(yz))+((([xu]z)y)[v])-(([xu](z[v]))y).
\end{eqnarray*}
(ii) Since
\begin{eqnarray*}
&&((((([xu]y)[v_1])[v_2]) \cdots [v_n])z)\\
&=&((((([xu]y)[v_1]) \cdots [v_{n-1}])z)[v_n])-(((([xu]y)[v_1])
\cdots [v_{n-1}])(z[v_n]))
\end{eqnarray*}
and
\begin{eqnarray*}
&& ([xu](((y[v_1]) \cdots [v_n])z))\\
&=&([xu](((y[v_1]) \cdots [v_{n-1}])z)[v_n])-([xu]((y[v_1]) \cdots
[v_{n-1}])(z[v_n]))\\
&=&(([xu](((y[v_1]) \cdots
[v_{n-1}])z))[v_n])-(([xu][v_n])(((y[v_1]) \cdots
[v_{n-1}])z))\\
&&-(([xu]((y[v_1]) \cdots
[v_{n-1}]))(z[v_n]))+(([xu](z[v_n]))((y[v_1]) \cdots [v_{n-1}]))\\
&=&(([xu](((y[v_1]) \cdots
[v_{n-1}])z))[v_n])-(([xu][v_n])(((y[v_1]) \cdots
[v_{n-1}])z)\\
&&-((([xu]((y[v_1]) \cdots [v_{n-2}]))[v_{n-1}])(z[v_n]))
+((([xu][v_{n-1}])((y[v_1]) \cdots [v_{n-2}]))(z[v_n]))\\
&&+(([xu](z[v_n]))((y[v_1]) \cdots [v_{n-1}]))\\
&=&(([xu](((y[v_1]) \cdots
[v_{n-1}])z))[v_n])-(([xu][v_n])(((y[v_1]) \cdots
[v_{n-1}])z)\\
&&-(((([xu]((y[v_1]) \cdots
[v_{n-3}]))[v_{n-2}])[v_{n-1}])(z[v_n]))\\
&&+(((([xu][v_{n-2}])((y[v_1]) \cdots [v_{n-3}]))[v_{n-1}])(z[v_n]))\\
&&+((([xu][v_{n-1}])((y[v_1]) \cdots [v_{n-2}]))(z[v_n]))
+(([xu](z[v_n]))((y[v_1]) \cdots [v_{n-1}]))\\
&=&\cdots\cdots \\
&=&(([xu](((y[v_1]) \cdots
[v_{n-1}])z))[v_n])-(([xu][v_n])(((y[v_1]) \cdots
[v_{n-1}])z)\\
&&-(((([xu]y)[v_1]) \cdots
[v_{n-1}])(z[v_n]))+(([xu](z[v_n]))((y[v_1]) \cdots [v_{n-1}]))\\
&&+\sum_{i=1}^{n-1}((((([xu][v_i])((((y[v_1])[v_2]) \cdots
[v_{i-1}])))[v_{i+1}])\cdots[v_{n-1}])(z[v_n])),
\end{eqnarray*}
we can get the result. \hfill $\blacksquare$

\begin{theorem}
Let $Lie(X)$ be the free Lie algebra generated by $X$ over a
commutative ring $K$ with unit. Then with deg-lex ordering on $X^*$,
the set $S=\{[xuy] \ | \ x,y,z\in X,\ u\in X^*, x\rhd y\rhd
supp(u)\}$ forms a Gr\"{o}bner-Shirshov basis in $Lie(X)$. As a
result, $Irr(S)=\{[u] \ | \ [u] \mbox{ is a } NLSW, \
u\neq{a\bar{s}b}, \ s\in{S},\ a,b\in{X^*}\}$ is a $K$-basis of the
free partially commutative Lie algebra $Lie(X| \vartheta)=Lie(X|S)$.
\end{theorem}
\textbf{Proof}\ Let us check all the possible compositions. The
ambiguities $w$ of  all possible compositions are:
\begin{enumerate}

\item[(i)] \ $w=xuzvy,\ x,y,z\in X,\  u,v \in X^*,\ x\rhd y \rhd supp(uzv),\ x\rhd z \rhd supp(u)$.

\item[(ii)] \ $w=xu|_{x_1u_1y_1}y,\ x,y,x_1,y_1\in X,\ u,u_1 \in X^*,\ x\rhd y \rhd supp(u),\ x_1\rhd y_1
\rhd supp(u_1)$.

\item[(iii)] \ $w=xuyvz,\ x,y,z\in X,\  u,v \in X^*,\ x\rhd y \rhd supp(u),\  y\rhd z \rhd supp(v)$.

\end{enumerate}
Now we prove that all the compositions are trivial.

For (i),  let $f=[xuzvy]$, $g=[xuz]$, $x,y,z\in X,\  u,v \in X^*,\
x\rhd y \rhd supp(uzv),\ x\rhd z \rhd supp(u)$. Suppose $v=v_1v_2
\cdots v_n$, where $v_i$ is an $ALSW$, \ $1 \leq i \leq n$, \ \ and
\ $v_1 \leq v_2 \leq \cdots \leq v_n$. Then $w=xuzvy$ and $
(f,g)_w=[xuzvy]-[xuzvy]_{xuz}=([xuzv]y)-(((([xuz][v_1])[v_2]) \cdots
[v_n])y) =([xuzv]-((([xuz][v_1])[v_2]) \cdots [v_n])y). $ By Lemmas
\ref{l3.5} and  \ref{l3.6}, we have that
$[xuzv]-((([xuz][v_1])[v_2]) \cdots [v_n])=\sum_{i=1}^m
\alpha_i[xw_i]$ where $\alpha_i\in k$,  $xw_i$ is an $ALSW$,
$supp(xw_i)=supp(xuzv)$, $|xw_i|_X=|xuzv|_X$  and
$\overline{[xw_i]}< xuzv$ for any $i=1,2,\ldots m$. Then
$$
(f,g)_w=((\sum_{i=1}^m\alpha_i [xw_i])y)= \sum_{i=1}^m\alpha_i
([xw_i]y)=\sum_{i=1}^m\alpha_i [xw_iy]\equiv0\ \  mod(S, w).
$$

For (ii),  let $f=[xu|_{x_1u_1y_1}y]$, $g=[x_1u_1y_1]$,
$x,y,x_1,y_1\in X,\ u,u_1 \in X^*,\ x\rhd y \rhd supp(u),\ x_1\rhd
y_1 \rhd supp(u_1)$. Then $w=xu|_{x_1u_1y_1}y$ and
$(f,g)_w=[xu|_{x_1u_1y_1}y]-[xu|_{x_1u_1y_1}y]_{x_1u_1y_1}=
([xu|_{x_1u_1y_1}]y)-(([xu|_{x_1u_1y_1}]_{x_1u_1y_1})y)=
(([xu|_{x_1u_1y_1}]-([xu|_{x_1u_1y_1}]_{x_1u_1y_1}))y)$.  By Lemmas
\ref{l3.5} and \ref{l3.6}, we have that
$[xu|_{x_1u_1y_1}]-([xu|_{x_1u_1y_1}]_{x_1u_1y_1})=\sum_{i=1}^l\alpha_i
[xw_i]$ where $\alpha_i\in k$, $xw_i$ is an $ALSW$,
$supp(xw_i)=supp(xu|_{x_1u_1y_1})$,  $|xw_i|_X=|xu|_{x_1u_1y_1}|_X$
and $\overline{[xw_i]}< xu|_{x_1u_1y_1}$ for any $i=1,2,\ldots l$.
Then
$$
(f,g)_w=((\sum_{i=1}^l \alpha_i[xw_i])y)=\sum_{i=1}^l
\alpha_i([xw_i]y)=\sum_{i=1}^l \alpha_i[xw_iy]\equiv0\ \  mod(S, w).
$$

For (iii),  let $f=[xuy]$, $g=[yvz]$, $x,y,z\in X,\ u,v \in X^*,\
x\rhd y \rhd supp(u),\ y\rhd z \rhd supp(v)$. There are two cases to
consider.

(a)  $|v|=0$. Then $w=xuyz$ and $ (f,g)_w=[xuyz]_{xuy}-[xuyz]_{yz}=
(([xu]y)z)-([xu](yz))= (([xu]y)z)-(([xu]y)z)+(([xu]z)y)=(([xu]z)y).
$ By Lemmas \ref{l3.5} and  \ref{l3.6}, we have that
$([xu]z)=\sum_{i=1}^l \alpha_i[xw_i]$, where $\alpha_i\in k$, $xw_i$
is an $ALSW$, $supp(xw_i)=supp(xuz)$, $|xw_i|_X=|xuz|_X$ and
$\overline{[xw_i]}\leq \overline{([xu]z)}=xuz$ for any $i=1,2,\ldots
l$. Then
$$
(f,g)_w=((\sum_{i=1}^l \alpha_i[xw_i])y)=\sum_{i=1}^l
\alpha_i([xw_i]y)=\sum_{i=1}^l \alpha_i[xw_iy]\equiv0\ \  mod(S, w).
$$

(b)  $|v|\geq 1$. Suppose  $v=v_1v_2 \cdots v_n$, where $v_i$ is an
$ALSW$, \ $1 \leq i \leq n$, \ \ and \ $v_1 \leq v_2 \leq \cdots
\leq v_n$. Then $w=xuyvz$ and by Lemma \ref{5.14}, we have that

\begin{eqnarray*}
&&(f,g)_w=[xuyvz]_{xuy}-[xuyvz]_{yvz}\\
&=&((((([xu]y)[v_1])[v_2]) \cdots [v_n])z)-([xu]((((y[v_1])[v_2])
\cdots [v_n])z))\\
&=&(((((([xu]y)[v_1])[v_2]) \cdots
[v_{n-1}])z)[v_n])-(([xu]((((y[v_1])[v_2]) \cdots
[v_{n-1}])z))[v_n])\\
&&+(([xu][v_n])((((y[v_1])[v_2]) \cdots
[v_{n-1}])z))-(([xu](z[v_n]))(((y[v_1])[v_2])
\cdots [v_{n-1}]))\\
&&-\sum_{i=1}^{n-1}((((([xu][v_i])((((y[v_1])[v_2]) \cdots
[v_{i-1}])))[v_{i+1}])\cdots[v_{n-1}])(z[v_n]))\\
&=&(((((([xu]y)[v_1]) \cdots
[v_{n-2}])z)[v_{n-1}])[v_n])-((([xu](((y[v_1]) \cdots
[v_{n-2}])z))[v_{n-1}])[v_n])\\
&&+((([xu][v_{n-1}])(((y[v_1]) \cdots
[v_{n-2}])z))[v_n])-((([xu](z[v_{n-1}]))((y[v_1])
\cdots [v_{n-2}]))[v_n])\\
&&-\sum_{i=1}^{n-2}(((((([xu][v_i])(((y[v_1]) \cdots
[v_{i-1}])))[v_{i+1}])\cdots[v_{n-2}])(z[v_{n-1}]))[v_n])\\
&&+(([xu][v_n])(((y[v_1]) \cdots
[v_{n-1}])z))-(([xu](z[v_n]))((y[v_1])
\cdots [v_{n-1}]))\\
&&-\sum_{i=1}^{n-1}((((([xu][v_i])(((y[v_1]) \cdots
[v_{i-1}])))[v_{i+1}])\cdots[v_{n-1}])(z[v_n]))\\
&=&\cdots\cdots\\
&=&(((([xuy][v_1])z)[v_2]) \cdots [v_n])-((([xu]((y[v_1])z))[v_2])
\cdots [v_n])\\
&&+(((([xu][v_{2}])((y[v_1])z))[v_3])\cdots[v_n])-(((([xu](z[v_{2}]))(y[v_1]))[v_3])\cdots[v_n])\\
&&-((((([xu][v_1])y)(z[v_2]))
[v_{3}])\cdots[v_n])+\cdots\cdots\\
&&+((([xu][v_{n-1}])(((y[v_1]) \cdots
[v_{n-2}])z))[v_n])-((([xu](z[v_{n-1}]))((y[v_1])
\cdots [v_{n-2}]))[v_n])\\
&&-\sum_{i=1}^{n-2}(((((([xu][v_i])(((y[v_1]) \cdots
[v_{i-1}])))[v_{i+1}])\cdots[v_{n-2}])(z[v_{n-1}]))[v_n])\\
&&+(([xu][v_n])(((y[v_1]) \cdots
[v_{n-1}])z))-(([xu](z[v_n]))((y[v_1])
\cdots [v_{n-1}]))\\
&&-\sum_{i=1}^{n-1}((((([xu][v_i])(((y[v_1]) \cdots
[v_{i-1}])))[v_{i+1}])\cdots[v_{n-1}])(z[v_n]))\\
&=&(((([xu][v_1])(yz))[v_2])\cdots[v_n])+((((([xu]z)y)[v_1])[v_2])\cdots[v_n])\\
&&-(((([xu](z[v_1]))y)[v_2])\cdots[v_n])+(((([xu][v_{2}])((y[v_1])z))[v_3])\cdots[v_n])\\
&&-(((([xu](z[v_{2}]))(y[v_1]))[v_3])\cdots[v_n])-((((([xu][v_1])y)(z[v_2]))
[v_{3}])\cdots[v_n])\\
&&+\cdots\cdots+\\
&&+((([xu][v_{n-1}])(((y[v_1]) \cdots
[v_{n-2}])z))[v_n])\\
&&-((([xu](z[v_{n-1}]))((y[v_1])
\cdots [v_{n-2}]))[v_n])\\
&&-\sum_{i=1}^{n-2}(((((([xu][v_i])(((y[v_1]) \cdots
[v_{i-1}])))[v_{i+1}])\cdots[v_{n-2}])(z[v_{n-1}]))[v_n])\\
&&+(([xu][v_n])(((y[v_1]) \cdots
[v_{n-1}])z))-(([xu](z[v_n]))((y[v_1])
\cdots [v_{n-1}]))\\
&&-\sum_{i=1}^{n-1}((((([xu][v_i])(((y[v_1]) \cdots
[v_{i-1}])))[v_{i+1}])\cdots[v_{n-1}])(z[v_n]))\\
&\equiv&0\ \ mod(S, w).
\end{eqnarray*}

So, all compositions in $S$ are trivial.

Clearly, $Lie(X| \vartheta)=Lie(X| S)$.

The latter conclusion follows from the Composition-Diamond lemma for
Lie algebras (Lemma \ref{cdL}).  \hfill $\blacksquare$

\ \

\noindent{\bf Acknowledgement}: The authors would like to thank
Professor L.A. Bokut for his guidance, useful discussions and
enthusiastic encouragement in writing up this paper.

\end{document}